\documentclass{amsart}
\usepackage{amssymb}

\usepackage{graphicx}
\usepackage{epstopdf}

\input tcilatex

\begin{document}
\title{On a coalescence process and its branching genealogy}
\author{Nicolas Grosjean, Thierry Huillet}
\address{Laboratoire de Physique Th\'{e}orique et Mod\'{e}lisation \\
CNRS-UMR 8089 et Universit\'{e} de Cergy-Pontoise,\\
2 Avenue Adolphe Chauvin, 95302, Cergy-Pontoise, FRANCE\\
E-mail: Nicolas.Grosjean@u-cergy.fr, Thierry.Huillet@u-cergy.fr}
\maketitle

\begin{abstract}
We define and analyze a coalescent process as a recursive box-filling
process whose genealogy is given by an ancestral time-reversed,
time-inhomogeneous Bienyam\'{e}-Galton-Watson process. Special interest is
on the expected size of a typical box and its probability of being empty.
Special cases leading to exact asymptotic computations are investigated when
the coalescing mechanisms are either linear-fractional or quadratic.\newline

\end{abstract}

\section{Introduction}

We are interested in a model of particles iteratively coalescing at random
at different sites indexing (say) energy levels and possibly leading to
unoccupied sites free of particles and with no immigration of particles from
outside. This can also be described as a random reallocation of balls in
boxes possibly generating empty boxes, or as an iterative coarse-graining of
events in a renewal process possibly generating moments free of events, or
as a random walk of possibly merging particles. In such processes, there is
a competition between a random force enhancing the merging of particles at
some sites and a balancing random force whose effect is to generate empty
sites and this trade-off is controlled by a sequence of inhomogeneous
coalescing mechanisms. There is a strong analogy with the time-inhomogeneous
Bienaym\'{e}-Galton-Watson (BGW) branching processes but the coalescing
process is of a different nature. In fact, the genealogy of such coalescing
processes after a fixed number of steps, say $N$, is a time-reversed BGW
branching process conditioned on its number of offspring at the terminal
value $N$. For some remarkable coalescing mechanism sequences, the expected
occupancies of boxes together with the probability of a box being empty can
be computed and estimated in the large $N$ limit. This concerns the
linear-fractional (or $\theta $-coalescing) and the quadratic coalescing
mechanisms that we shall investigate in some detail. Connections, different
from ours, between coalescence Markov processes and the classical branching
theory (for various branching mechanisms) appeared in \cite{A}, \cite{A2}, 
\cite{Lam} and \cite{Sheth}, which deal with the probability that randomly
sampled individuals of a time-homogeneous BGW process, alive at some
generation, merge in the previous generation.

\section{A coalescence process and its branching genealogy}

\subsection{A coalescence process as a recursive box-filling process}

Let $\Bbb{N}_{0}:=\left\{ 0,1,2,..\right\} $ and $\Bbb{N}:=\left\{
1,2,..\right\} $. We introduce a Markovian coalescence process with
state-space $\Bbb{N}_{0}^{\Bbb{N}}$ defined as follows: consider an infinite
sequence $\left\{ K_{n}^{*}\left( i\right) ,\text{ }i\geq 1\right\} $ giving
the state of the process at step $n\in \Bbb{N}_{0}$; $K_{n}^{*}\left(
i\right) \in \Bbb{N}_{0}$ will be the number of balls at step $n$ in box
number $i$ (or at site $i$), in a scenario with infinitely many boxes
(sites). Suppose at step $n=0$, $\left\{ K_{0}^{*}\left( i\right) =1,i\geq
1\right\} $ (all boxes are filled just with one ball). For each $n\geq 1$,
let $\left\{ M_{n}\left( i\right) ,\text{ }i\geq 1\right\} $ be an infinite
sequence of independent and identically distributed (iid) $\Bbb{N}_{0}$%
-valued random variables indexed by $n$, and let $f_{n}\left( z\right) :=%
\mathbf{E}\left( z^{M_{n}}\right) $ be their common probability generating
function (pgf) for which it is assumed $f_{n}\left( 1\right) =1$ and $%
1>f_{n}\left( 0\right) >0$ for each $n\geq 1$. We shall call $f_{n}\left(
z\right) $ the coalescing mechanism at step $n$. Suppose $\left\{
M_{n}\left( i\right) ,\text{ }i\geq 1\right\} $ are also mutually
independent across $n$.

\begin{definition}
The updating mechanism of $\left\{ K_{n}^{*}\left( i\right) ,i\geq 1\right\} 
$ is defined as follows: 
\begin{equation}
K_{n+1}^{*}\left( i\right) =\sum_{m=1}^{M_{n+1}\left( i\right)
}K_{n}^{*}\left( m+\sum_{j=1}^{i-1}M_{n+1}\left( j\right) \right) \text{, }%
i\geq 1  \label{1.a}
\end{equation}
with the convention $K_{n+1}^{*}\left( i\right) =0$ if $M_{n+1}\left(
i\right) =0$.
\end{definition}

Note that a given box will never be empty for ever. We will find it also
useful below to consider this dynamics stopped up to some terminal value,
say $N$, of $n$. Let us comment the Markovian reallocation dynamics (\ref
{1.a}):\newline

Suppose $\left\{ K_{n}^{*}\left( i\right) ,\text{ }i\geq 1\right\} $ are
given. To compute the number of balls in box number $i$ at step $n+1$, draw
a sequence $\left\{ M_{n+1}\left( i\right) ,\text{ }i\geq 1\right\} $. Fill
the first box at step $n+1$ by coalescing (or merging) the balls in the $%
M_{n+1}\left( 1\right) $ first boxes of the $n$-th step, with the convention
that if $M_{n+1}\left( 1\right) =0$, an empty box is created instead. This
gives the size of the first box $K_{n+1}^{*}\left( 1\right) $ at step $n+1$.
To compute $K_{n+1}^{*}\left( 2\right) $, use the remaining boxes at step $n$
with $i>M_{n+1}\left( 1\right) $ and merge the $M_{n+1}\left( 2\right) $
first remaining boxes. And iterate the process to form $\left\{
K_{n+1}^{*}\left( i\right) ,\text{ }i\geq 1\right\} $. From (\ref{1.a}), it
is finally clear that balls filling up some box at some step all come up
from balls in some box at the preceding step: in our model of coalescence,
there is no incoming balls from outside and therefore the system is closed
or conservative; the adjunction of immigrants possibly filling up the empty
sites at each step looks a promising issue.

\begin{proposition}
The $\left\{ K_{n}^{*}\left( i\right) ,\text{ }i\geq 1\right\} $ are
mutually iid for each $n$. With, say, $K_{n}^{*}\overset{\text{law}}{=}%
K_{n}^{*}\left( 1\right) $ and $\left\{ K_{n}^{*\left( m\right) },\text{ }%
m\geq 1\right\} $ iid copies of $K_{n}^{*},$%
\begin{equation}
K_{n+1}^{*}\overset{\text{law}}{=}\sum_{m=1}^{M_{n+1}}K_{n}^{*\left(
m\right) }.  \label{1.b}
\end{equation}
\end{proposition}

\textbf{Proof:} Due to the mutual independence $\left\{ M_{n+1}\left(
i\right) ,\text{ }i\geq 1\right\} $ both in $i$ for each $n$ and then also
across $n$ and because a given box at step $n$ will contribute to the
construction of only one box at step $n+1$, by induction, the $\left\{
K_{n}^{*}\left( i\right) ,\text{ }i\geq 1\right\} $ are all iid for each $n$%
. And (\ref{1.b}) follows from (\ref{1.a}). Let then $\phi _{n}^{*}\left(
z\right) =\mathbf{E}\left( z^{K_{n}^{*}}\right) $ be the common typical pgf
of the $\left\{ K_{n}^{*}\left( i\right) ,\text{ }i\geq 1\right\} $. Then 
\begin{equation}
\left\{ 
\begin{array}{c}
\phi _{n+1}^{*}\left( z\right) =f_{n+1}\left( \phi _{n}^{*}\left( z\right)
\right) \text{, }n\geq 0\text{ or} \\ 
\phi _{n}^{*}\left( z\right) =f_{n}\left( f_{n-1}\left( ...f_{1}\left(
z\right) \right) \right) \text{, }n\geq 1
\end{array}
\right. ,  \label{1.0}
\end{equation}
gives the law of $K_{n}^{*}$ as $\mathbf{P}\left( K_{n}^{*}=j\right) =\left[
z^{j}\right] \phi _{n}^{*}\left( z\right) $, $j\geq 0$. $\square $

Of particular interest is the probability $\mathbf{P}\left(
K_{n}^{*}=0\right) =\phi _{n}^{*}\left( 0\right) $ to have a typical empty
box at step $n$, together with the mean size $\mathbf{E}\left(
K_{n}^{*}\right) $ of the typical box at step $n$. With $n_{1}>n_{2}$, we
can also define 
\begin{equation*}
\Phi _{n_{1},n_{2}}^{*}\left( z\right) :=f_{n_{1}}\left( f_{n_{1}-1}\left(
...f_{n_{2}+1}\left( z\right) \right) \right) .
\end{equation*}
Then, with $n_{1}>n_{2}>n_{3}$, $\Phi _{n_{1},n_{3}}^{*}\left( z\right)
=\Phi _{n_{1},n_{2}}^{*}\left( z\right) \circ \Phi _{n_{2},n_{3}}^{*}\left(
z\right) :=\Phi _{n_{1},n_{2}}^{*}\left( \Phi _{n_{2},n_{3}}^{*}\left(
z\right) \right) $, a concatenation property of the time-inhomogeneous
coalescent process $\left\{ K_{n}^{*}\right\} $. It holds $\Phi
_{n,0}^{*}\left( z\right) =\phi _{n}^{*}\left( z\right) $.

\begin{figure}[tbp]
\label{lafig1}
\par
\begin{center}
\includegraphics[scale=.35,clip]{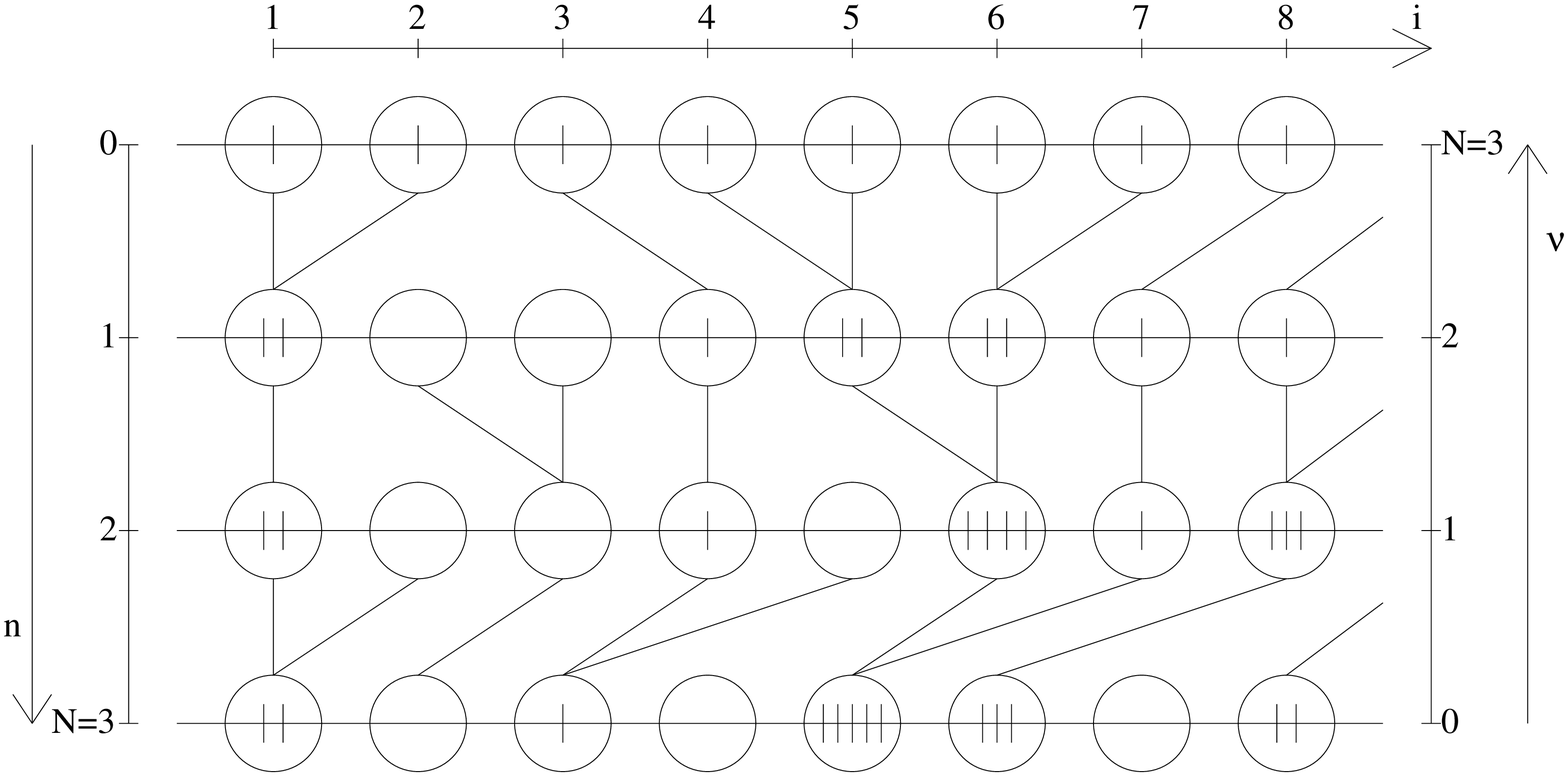}
\end{center}
\caption{Example of a binary coalescence process with the following
realizations of $M_{n}(i)$: \newline
$M_{1}(1)=2$, $M_{1}(2)=0$, $M_{1}(3)=0$, $M_{1}(4)=1$, $M_{1}(5)=2$, $%
M_{1}(6)=2$, $M_{1}(7)=1$, $M_{1}(8)=1$,...; \newline
$M_{2}(1)=1$, $M_{2}(2)=0$, $M_{2}(3)=2$, $M_{2}(4)=1$, $M_{2}(5)=0$, $%
M_{2}(6)=2$, $M_{2}(7)=1$, $M_{2}(8)=2$,...; \newline
$M_{3}(1)=2$, $M_{3}(2)=1$, $M_{3}(3)=2$, $M_{3}(4)=0$, $M_{3}(5)=2$, $%
M_{3}(6)=1$, $M_{3}(7)=0$, $M_{3}(8)=1$,.... }
\end{figure}

Although the $\left\{ K_{n}^{*}\left( i\right) ,\text{ }i\geq 1\right\} $
are all iid for each $n$, for any $\left( i,j\right) $, $K_{n}^{*}\left(
i\right) $ and $K_{n+1}^{*}\left( j\right) $ are not necessarily
independent: for instance if $K_{n}^{*}\left( 1\right) $ is large, $%
K_{n+1}^{*}\left( 1\right) $ is more likely to be large.

We also observe that for any finite $n$, $\phi _{n}^{*}\left( 0\right) $
cannot be equal to $1$. Indeed, if this were the case, there would exist $%
n_{*}\leq n$ such that $\phi _{n_{*}}^{*}\left( 0\right) =1$ and $\phi
_{n_{*}-1}^{*}\left( 0\right) <1$ (recalling $\phi _{0}^{*}\left( 0\right)
=0 $). Now $\phi _{n_{*}}^{*}\left( 0\right) =f_{n_{*}}\left( \phi
_{n_{*}-1}^{*}\left( 0\right) \right) $ and $f_{n_{*}}$ is invertible
because $f_{n_{*}}$ is not constant equal to $1$; therefore $\phi
_{n_{*}-1}^{*}\left( 0\right) =1$ which contradicts the existence of $n_{*}$%
: for any $n$, all $\left\{ K_{n}^{*}\left( i\right) ,\text{ }i\geq
1\right\} $ cannot be simultaneously equal to $0$. Moreover, the sequence $%
\phi _{n}^{*}\left( 0\right) $ is neither necessarily increasing nor
decreasing with $n$ (see the examples below).

By independence of the $K_{n}^{*}\left( i\right) $s, in view of $0<\phi
_{n}^{*}\left( 0\right) <1$, for each fixed $n\in \left\{ 1,...,N\right\} $,
we also have

\begin{equation*}
\mathbf{P}\left( \min_{i\geq 1}K_{n}^{*}\left( i\right) >0\right) =\left(
1-\phi _{n}^{*}\left( 0\right) \right) ^{\infty }=0,
\end{equation*}
(at least one box is empty) and by strong law of large numbers, $\mathbf{P}%
\left( K_{n}^{*}=j\right) =\lim_{I\rightarrow \infty }\frac{1}{I}%
\sum_{i=1}^{I}\mathbf{1}_{\left\{ K_{n}^{*}\left( i\right) =j\right\} }$,
for any $j\geq 0$, together with $\mathbf{E}\left( K_{n}^{*}\right)
=\lim_{I\rightarrow \infty }\frac{1}{I}\sum_{i=1}^{I}K_{n}^{*}\left(
i\right) $, if this quantity exists (is finite). In particular, $\mathbf{P}%
\left( K_{n}^{*}=0\right) =\phi _{n}^{*}\left( 0\right) $. Note finally, 
\begin{equation*}
\lim_{I\rightarrow \infty }\frac{1}{I}\sum_{i=1}^{I}M_{n+1}\left( i\right) =%
\mathbf{E}\left( M_{n+1}\right) =f_{n+1}^{\prime }\left( 1\right) ,
\end{equation*}
and $\mathbf{E}\left( K_{N}^{*}\right) =\phi _{N}^{*\prime }\left( 1\right)
=\prod_{n=1}^{N}f_{n}^{\prime }\left( 1\right) =\prod_{n=1}^{N}\mathbf{E}%
\left( M_{n}\right) $, if these quantities exist (are finite).

\begin{definition}
We shall say that we have a subcritical, critical or supercritical
coalescence process if the (geometric mean) limit $\mu =\lim_{N\rightarrow
\infty }\mathbf{E}\left( K_{N}^{*}\right) ^{1/N}$ exists and $\mu <1$, $=1$, 
$>1$, respectively. If for each $n$, $\mathbf{E}\left( M_{n}\right) =\infty $
so that $\mathbf{E}\left( K_{N}^{*}\right) =\infty $, we shall put $\mu
:=\infty $ and the process will be said strongly supercritical.
\end{definition}

In the (sub-)critical cases, we expect $\mathbf{P}\left( K_{N}^{*}=0\right) 
\underset{N\rightarrow \infty }{\rightarrow }1$, and of interest is the rate
at which $\mathbf{P}\left( K_{N}^{*}>0\right) $ tends to $0$. In the
supercritical case, we expect $\mathbf{P}\left( K_{N}^{*}=0\right) \underset{%
N\rightarrow \infty }{\rightarrow }\rho \in \left[ 0,1\right) $ and of
interest is when $\rho >0$ and the rates at which $\mathbf{P}\left(
K_{N}^{*}>0\right) $ tends to $1-\rho \geq 0$. In this setup, we will give
below one example for which $\rho \in \left( 0,1\right) $ (the supercritical 
$\theta $-coalescing model below) and another one for which $\rho =0$ (the
supercritical quadratic branching model).\newline

In Figure $1$ we give a partial view of a sample realization of the binary
coalescence process. Picking the size $j=5$ node at terminal step $N=3$ and
moving backwards the steps, we see a tree that goes up to step $0$ with $5$
leaves of size $1$. Picking the first size $j=0$ node at step $N=3$ and
moving backwards, we see a tree made of empty nodes only that goes extinct
at step $1$. The ancestral tree of the first size $j=1$ node at step $N=3$
has an empty node at step $2$. It goes up to the step $0$ with $1$ leaf of
size $1$.

\begin{remark}
\textbf{\ }$\left( i\right) $\emph{\ If the law of }$M_{n+1}\left( i\right) $%
\emph{\ is the same for each }$n\geq 0$\emph{\ and }$i\geq 1$\emph{\ (with
the obvious simplifications of the construction), we shall speak of a
homogeneous coalescence process with homogeneous coalescing mechanism }$f$%
\emph{, independent of }$n$\emph{.\newline
}

$\left( ii\right) $\emph{\ A slight extension of (\ref{1.a}): let }$\left(
t\left( i\right) ,i\geq 1\right) $\emph{\ be a positive sequence, possibly
random and iid, representing say time lags (energy spacings or gaps). Let }$%
s\left( i\right) =\sum_{j=1}^{i}t\left( j\right) $\emph{, }$i\geq 1$\emph{\
be the moments of events (the energy levels). Defining }$\left\{
K_{n}^{*}\left( s\left( i\right) \right) ,i\geq 1\right\} $\emph{\ to be the
number of events (or particles) at times (or energy levels) }$s\left(
i\right) $\emph{, }$\left\{ K_{0}^{*}\left( s\left( i\right) \right)
=1,i\geq 1\right\} $\emph{\ is a renewal sequence of events and } 
\begin{equation*}
K_{n+1}^{*}\left( s\left( i\right) \right) =\sum_{m=1}^{M_{n+1}\left(
i\right) }K_{n}^{*}\left( s\left( m+\sum_{j=1}^{i-1}M_{n+1}\left( j\right)
\right) \right) \text{, }i\geq 1
\end{equation*}
\emph{gives the way the renewal sequence is updated by iterated coalescence
of events in the previous step, as a coarse-graining process possibly
generating empty sites. As in (\ref{1.a}), this process is driven by the
coalescing mechanism sequences.}
\end{remark}

\subsection{Genealogy as a time-reversed inhomogeneous BGW branching process}

Consider a discrete-time and time-inhomogeneous BGW branching process \cite
{AN}, \cite{Harris}, \cite{Jag}, \cite{tu}, whose reproduction laws are
given by the probability systems $\mathbf{P}\left( M_{n}=m\right) =\pi
_{n}\left( m\right) $, $m,n\geq 1$ for the number $M_{n}$ of offspring per
capita produced at generation $n$ to form generation $n+1$. We assume $\pi
_{n}\left( 0\right) >0$ so that the process can go extinct. We let $%
f_{n}\left( z\right) =\mathbf{E}\left( z^{M_{n}}\right) =\sum_{m\geq 0}\pi
_{n}\left( m\right) z^{m}$ be the probability generating function of $M_{n}$
(the branching mechanism) and we assume $f_{n}\left( 1\right) =1$ for each $%
n\geq 1$. With \textrm{K}$_{n}$ the number of individuals alive at
generation $n$ given \textrm{K}$_{0}=1$, we have \textrm{K}$%
_{n+1}=\sum_{k=1}^{\mathrm{K}_{n}}M_{n+1}^{\left( k\right) }$ where $%
M_{n+1}^{\left( k\right) }$ are iid copies of $M_{n+1}$ giving the offspring
number of each individual of the population alive at generation $n$. Thus, 
\begin{equation}
\phi _{n+1}\left( z\right) =\phi _{n}\left( f_{n+1}\left( z\right) \right)
,\phi _{0}\left( z\right) =z,  \label{1.1a}
\end{equation}
leading classically to 
\begin{equation*}
\mathbf{E}\left( z^{\mathrm{K}_{n}}\right) :=\phi _{n}\left( z\right)
=f_{1}\left( f_{2}\left( ...f_{n}\left( z\right) \right) \right) ,n\geq 1.
\end{equation*}
$\phi _{n}\left( z\right) $ is the $n$-th composition of the different $%
f_{m}\left( z\right) $ in the reverse order to the one giving $\phi _{n}^{*}$
in (\ref{1.0}). Recall that with $n_{1}<n_{2}<n_{3}$, if 
\begin{equation*}
\Phi _{n_{1},n_{2}}\left( z\right) :=f_{n_{1}+1}\left( f_{n_{1}+2}\left(
...f_{n_{2}}\left( z\right) \right) \right) ,
\end{equation*}
then $\Phi _{n_{1},n_{3}}\left( z\right) =\Phi _{n_{1},n_{2}}\left( z\right)
\circ \Phi _{n_{2},n_{3}}\left( z\right) :=\Phi _{n_{1},n_{2}}\left( \Phi
_{n_{2},n_{3}}\left( z\right) \right) $, the standard concatenation property
of the time-inhomogeneous BGW process $\left\{ \mathrm{K}_{n}\right\} $.

When running the above branching process $\left\{ \mathrm{K}_{n}\right\} $
backward in time starting from a terminal generation say $n=N$ to $n=0$ for
any given $N>0$, we get a new (truncated) branching process say $\left\{
K_{\nu },\text{ }0\leq \nu \leq N\right\} $ defined as follows:

\begin{definition}
The process $\left\{ K_{\nu },\text{ }0\leq \nu \leq N\right\} $ is defined
by its pgf at generation number $\nu $: $\mathbf{E}\left( z^{K_{0}}\right)
=z $ and 
\begin{equation}
\text{ }\mathbf{E}\left( z^{K_{\nu }}\right) :=\Phi _{N,N-\nu }^{*}\left(
z\right) =f_{N}\left( f_{N-1}\left( ...f_{N-\nu +1}\left( z\right) \right)
\right) ,\text{ }1\leq \nu \leq N,  \label{1.1}
\end{equation}
so with reversed pgf sequence $g_{1}=f_{N},...,g_{\nu }=f_{N-\nu
+1},...,g_{N}=f_{1}$.
\end{definition}

This suggests that $\left\{ K_{n}^{*},\text{ }0\leq n\leq N\right\} $ and $%
\left\{ K_{\nu },\text{ }0\leq \nu \leq N\right\} $ could be time-reversed
of one another (in a different sense to the time-reversal defined in \cite
{Esty}), with the latter being the ancestral process of the former. Note in
this respect that $\nu $ in Figure $1$ runs from bottom to top of the graph,
starting from $\nu =0$.\newline

The relation between the cluster coalescence process $\left\{ K_{n}^{*},%
\text{ }0\leq n\leq N\right\} $ considered up to a given terminal generation 
$N$ and the BGW process $\left\{ K_{\nu },\text{ }0\leq \nu \leq N\right\} $
in (\ref{1.1}), as time-reversed processes, is indeed as follows:

\begin{theorem}
At given terminal step $N>1$ of the process $\left\{ K_{n}^{*}\left(
i\right) ,\text{ }i\geq 1\right\} $, pick any box number $i$ and suppose it
has size $j$, an event occurring with probability $\mathbf{P}\left(
K_{N}^{*}=j\right) $. Then the ancestral tree of this size $j$ box, running
backward in $n$ (so forward in $\nu $), will be the one of a BGW process as
in (\ref{1.1}), started from a single individual (the selected box) but
conditioned on having $j$ descendants at generation $\nu =N$.
\end{theorem}

\textbf{Proof: }$\left( i\right) $ The creation of an empty box at step $N$,
running forward in $n$, is understood as the extinction before time $N$ of
its ancestral tree as a time-reversed BGW process. This empty box has indeed
no non-empty offspring (parents) at generation $N-1$ although it can have
one or more empty boxes offspring. If the randomly chosen box is empty, its
ancestral tree is therefore made of empty nodes only and its height will be
strictly less than $N$. It is the one of a time-reversed BGW process started
from a single individual and conditioned on having no descendants
(extinction) at generation $\nu =N$, meaning $K_{N}=0$.

$\left( ii\right) $ The coalescence of $m\geq 1$ boxes into one single box
forward in $n$ at step $N$ is understood as this box giving birth to $m$
offspring at step $N-1$, backward in $n$ (forward in $\nu $), one of which
at least being non empty. If the randomly chosen box is non-empty (of size $%
j>0$) indeed, the content of this box (its size) is read from its number of
ancestors at generation $N$: its BGW ancestral tree has exactly $j$ leaves
at generation $\nu =N$ and it will have height $N$. We conclude for instance
that the law of the number $K_{\nu }$, $1\leq \nu \leq N$, of ancestors at
generation $\nu $ of a selected size $j$ box at step $N$ produced by the
coalescence process (\ref{1.a}), is given from (\ref{1.1}) by 
\begin{equation*}
\mathbf{P}\left( K_{\nu }=k\mid K_{N}=j\right) =\frac{\left[
z_{1}^{k}z_{2}^{j}\right] \Phi _{N,N-\nu }^{*}\left( z_{1}\Phi _{N-\nu
,0}^{*}\left( z_{2}\right) \right) }{\left[ z_{2}^{j}\right] \Phi
_{N,0}^{*}\left( z_{2}\right) },
\end{equation*}
recalling $\mathbf{P}\left( K_{\nu }=k,K_{N}=j\right) =\left[
z_{1}^{k}z_{2}^{j}\right] \Phi _{\nu ,0}^{*}\left( z_{1}\Phi _{N,\nu
}^{*}\left( z_{2}\right) \right) $. As required, we have 
\begin{equation*}
\sum_{j\geq 0}\mathbf{P}\left( K_{N}^{*}=j\right) \mathbf{E}\left(
z_{1}^{K_{\nu }}\mid K_{N}=j\right) =\sum_{j\geq 0}\left[ z_{2}^{j}\right]
\phi _{N}^{*}\left( z_{2}\right) \frac{\left[ z_{2}^{j}\right] \Phi
_{N,N-\nu }^{*}\left( z_{1}\Phi _{N-\nu ,0}^{*}\left( z_{2}\right) \right) }{%
\left[ z_{2}^{j}\right] \Phi _{N,0}^{*}\left( z_{2}\right) }
\end{equation*}
\begin{equation*}
=\sum_{j\geq 0}\left[ z_{2}^{j}\right] \Phi _{N,N-\nu }^{*}\left( z_{1}\Phi
_{N-\nu ,0}^{*}\left( z_{2}\right) \right) =\Phi _{N,N-\nu }^{*}\left(
z_{1}\Phi _{N-\nu ,0}^{*}\left( 1\right) \right)
\end{equation*}
\begin{equation*}
=\Phi _{N,N-\nu }^{*}\left( z_{1}\right) =f_{N}\left( f_{N-1}\left(
...f_{N-\nu +1}\left( z_{1}\right) \right) \right) =\mathbf{E}\left(
z_{1}^{K_{\nu }}\right) .\text{ }\square
\end{equation*}

\begin{remark}
\emph{Here is another connection of }$\left\{ K_{n}^{*},\text{ }0\leq n\leq
N\right\} $\emph{\ with the BGW process }$\left\{ K_{\nu },\text{ }0\leq \nu
\leq N\right\} $\emph{: in the civil system of kinship degree (see for
example \cite{www}), the degree of kinship between two individuals of the
same family is the number of edges between these individuals (an individual
is related to itself by a degree }$0$\emph{). The law of the number of
individuals of the generation }$N$\emph{\ that are related to some
individual at the same generation by a degree of kinship at most }$2n$\emph{%
\ is the same as the one of a coalescence process }$K_{n}^{*}$\emph{.}
\end{remark}

\section{Examples}

There are classes of discrete branching and/or coalescence processes for
which the pgf $\phi _{N}\left( z\right) $ of \textrm{K}$_{N}$ and/or the pgf 
$\phi _{N}^{*}\left( z\right) $ of $K_{N}^{*}$ is exactly computable,
thereby making the above computations concrete and somehow explicit, in the
large $N$ limit. One is related to the linear-fractional branching mechanism 
\cite{Sag} (so called here $\theta $-branching and that we shall rename $%
\theta $-coalescing in view of our time-reversal), the other one being the
quadratic coalescing mechanism which we shall define.

\subsection{The $\theta $-coalescing mechanism}

With $\left| \theta \right| \leq 1$, $a_{n},b_{n}>0,$ consider the $\theta $%
-coalescing mechanism model, defined by 
\begin{equation}
\left\{ 
\begin{array}{c}
f_{n}\left( z\right) =z_{c}-\left( a_{n}\left( z_{c}-z\right) ^{-\theta
}+b_{n}\right) ^{-1/\theta }\text{ or} \\ 
\left( z_{c}-f_{n}\left( z\right) \right) ^{-\theta }=a_{n}\left(
z_{c}-z\right) ^{-\theta }+b_{n},
\end{array}
\right. ,  \label{lnf}
\end{equation}
and for those values of $a_{n},b_{n}>0$ and $z_{c}\geq 1$, for which $f_{n}$
is a pgf with $f_{n}\left( 1\right) =1$. This is an inhomogeneous version of
a pgf family defined in \cite{Sag}.

\begin{proposition}
This family is stable under composition. Indeed, 
\begin{equation*}
\phi _{n}^{*}\left( z\right) =f_{n}\left( ...f_{1}\left( z\right) \right)
=z_{c}-\left( A_{n}^{*}\left( z_{c}-z\right) ^{-\theta }+B_{n}^{*}\right)
^{-1/\theta }
\end{equation*}
is a pgf with $A_{n}^{*}=\prod_{m=1}^{n}a_{m}$ and $%
B_{n}^{*}=B_{n}^{*}=b_{n}+...+b_{k}\prod_{l=k+1}^{n}a_{l}+...+b_{1}%
\prod_{l=2}^{n}a_{l}$. We also have 
\begin{equation}
\mathbf{P}\left( K_{N}^{*}=0\right) =\phi _{N}^{*}\left( 0\right)
=z_{c}-\left( A_{N}^{*}z_{c}^{-\theta }+B_{N}^{*}\right) ^{-1/\theta },
\label{eb1}
\end{equation}
and if the means exist, $f_{n}^{\prime }\left( 1\right) =a_{n}\left(
a_{n}+b_{n}\left( z_{c}-1\right) ^{\theta }\right) ^{-\left( 1+\theta
\right) /\theta }=a_{n}$ because $f_{n}\left( 1\right) =1$ entails $%
b_{n}=\left( z_{c}-1\right) ^{-\theta }\left( 1-a_{n}\right) $. Thus, 
\begin{equation}
\mathbf{E}\left( K_{N}^{*}\right) =\phi _{N}^{*\prime }\left( 1\right)
=\prod_{n=1}^{N}f_{N}^{\prime }\left( 1\right) =A_{N}^{*}\text{ with }%
\mathbf{E}\left( K_{N}^{*}\right) ^{1/n}=A_{N}^{*1/N}.  \label{mv1}
\end{equation}
\end{proposition}

For the sake of simplicity, we shall limit ourselves in the following
examples to the cases $z_{c}=1$ and $\theta \in \left\{ -1,0,1\right\} $.

\subsubsection{$\theta =-1$ \textbf{(affine case)}}

Here $f_{n}\left( z\right) =1-\left( a_{n}\left( 1-z\right) +b_{n}\right) $
and $f_{n}\left( 1\right) =1$ imposes $b_{n}=0$ and $a_{n}\in \left(
0,1\right) $. We have 
\begin{equation}
\phi _{n}^{*}\left( z\right) =1-A_{n}^{*}+A_{n}^{*}z.  \label{pgf1}
\end{equation}

- Take for example $a_{n}=\mu \left( 1-1/\left( n+1\right) ^{\alpha }\right) 
$ with $\alpha >0$.

If $\alpha =1$: $\mathbf{E}\left( K_{N}^{*}\right)
=\prod_{n=1}^{N}f_{n}^{\prime }\left( 1\right)
=\prod_{n=1}^{N}a_{n}=a^{N}/\left( n+1\right) $ so that $\mathbf{E}\left(
K_{N}^{*}\right) ^{1/N}\rightarrow \mu =a$. If $\mu =1$ ($\mu <1$) we get a
critical (subcritical) coalescence process. Here, $\mathbf{P}\left(
K_{N}^{*}=0\right) =\phi _{N}^{*}\left( 0\right)
=1-\prod_{n=1}^{N}a_{n}=1-\mu ^{N}/\left( N+1\right) \rightarrow 1.$
Asymptotically, the typical box gets empty with probability $1$, at rate $%
\mu ^{N}/\left( N+1\right) $ and $\mathbf{P}\left( K_{N}^{*}=0\right) $ is
increasing with $N$ for the affine model.

If $\alpha >1$: $\mathbf{E}\left( K_{N}^{*}\right) ^{1/N}\rightarrow \mu =a$%
. And $\mathbf{P}\left( K_{N}^{*}=0\right) =\phi _{N}^{*}\left( 0\right)
=1-\prod_{n=1}^{N}a_{n}=1-\mu ^{N}\frac{\zeta \left( \alpha \right) -1}{N}%
\rightarrow 1$, where $\zeta \left( \alpha \right) $ is the zeta function.

If $\alpha <1$: $\mathbf{P}\left( K_{N}^{*}=0\right) =\phi _{N}^{*}\left(
0\right) =1-\prod_{n=1}^{N}a_{n}=1-\frac{1}{\left( 1-\alpha \right) }%
N^{-\alpha }\mu ^{N}\rightarrow 1$. \newline

- More generally, take $a_{n}=a\left( 1-\varepsilon _{n}\right) $ where $%
\varepsilon _{n}\in \left[ 0,1\right) $ and $\varepsilon _{n}\underset{%
n\rightarrow \infty }{\rightarrow }0$. Then, it holds that the Ces\`{a}ro
mean $\left( \sum_{m=1}^{n}\varepsilon _{m}\right) /n\underset{n\rightarrow
\infty }{\rightarrow }0$. Thus $\mathbf{E}\left( K_{N}^{*}\right)
^{1/N}\rightarrow \mu =a$ and $\mathbf{P}\left( K_{N}^{*}=0\right)
\rightarrow 1$, the rates of which depending on the way $\varepsilon _{n}%
\underset{n\rightarrow \infty }{\rightarrow }0$.

\subsubsection{$\theta =1$ \textbf{(standard homographic case)}}

Here $f_{n}\left( z\right) $ takes the form 
\begin{equation}
f_{n}\left( z\right) =\left( 1-\frac{1}{a_{n}+b_{n}}\right) +\frac{1}{%
a_{n}+b_{n}}\frac{z}{1+b_{n}/a_{n}\left( 1-z\right) },  \label{bm1}
\end{equation}
where $a_{n},b_{n}>0$ and $a_{n}+b_{n}>1$. We have 
\begin{equation}
\phi _{n}^{*}\left( z\right) =1-\frac{1}{A_{n}^{*}+B_{n}^{*}}+\frac{1}{%
A_{n}^{*}+B_{n}^{*}}\frac{z}{1+B_{n}^{*}/A_{n}^{*}\left( 1-z\right) },
\label{pgf2}
\end{equation}
where $A_{n}^{*},B_{n}^{*}>0$ and $A_{n}^{*}+B_{n}^{*}>1$. We have $%
f_{n}^{\prime }\left( 1\right) =1/a_{n}$ and taking $a_{n}=a\left(
1-\varepsilon _{n}\right) $ where $\varepsilon _{n}\underset{n\rightarrow
\infty }{\rightarrow }0$, $\mathbf{E}\left( K_{N}^{*}\right)
=\prod_{n=1}^{N}f_{n}^{\prime }\left( 1\right) =\prod_{n=1}^{N}a_{n}^{-1}$
so that $\mathbf{E}\left( K_{N}^{*}\right) ^{1/N}\rightarrow \mu =1/a$.
Depending on $\mu <1$, $\mu =1$ or $\mu >1$, we get a subcritical, critical
or supercritical coalescence process. Moreover, 
\begin{equation}
\mathbf{P}\left( K_{N}^{*}=0\right) =\phi _{N}^{*}\left( 0\right) =1-\frac{1%
}{A_{N}^{*}+B_{N}^{*}}  \label{eb2}
\end{equation}
now depends on the choice of the sequence $b_{n}>1-a_{n}.$

Taking $b_{n}=b\left( 1-a_{n}\right) $ with $b>1$ if $a_{n}<1$ or $b<1$ if $%
a_{n}>1$, so that $a_{n}+b_{n}>1$ always, we get telescoping sums leading to 
$B_{N}^{*}=b\left( 1-A_{n}^{*}\right) $. Thus, 
\begin{equation*}
\mathbf{P}\left( K_{N}^{*}=0\right) =\frac{\left( b-1\right) \left(
1-A_{N}^{*}\right) }{A_{N}^{*}+b\left( 1-A_{N}^{*}\right) }<1\text{, }%
b\lessgtr 1\text{, if }A_{N}^{*}\gtrless 1.
\end{equation*}

If $a>1$ (a subcritical case because $\mu <1$), take for example $%
a_{n}=a\left( 1+\varepsilon _{n}\right) $ where $\varepsilon _{n}>0$ and $%
\varepsilon _{n}\underset{n\rightarrow \infty }{\rightarrow }0$ fast enough.
We get $\mathbf{P}\left( K_{N}^{*}=0\right) \rightarrow 1$ and with $%
A_{N}^{*-1/N}\rightarrow 1/a$: $\mathbf{P}\left( K_{N}^{*}>0\right) \sim 
\frac{1}{1-b}A_{N}^{*-1}\rightarrow 0.$

If $a<1$ (a supercritical case because $\mu >1$), take for example $%
a_{n}=a\left( 1-\varepsilon _{n}\right) $ where $\varepsilon _{n}\in \left(
0,1\right) $ and $\varepsilon _{n}\underset{n\rightarrow \infty }{%
\rightarrow }0$ fast enough. We get $\mathbf{P}\left( K_{N}^{*}=0\right)
\rightarrow \rho =\left( b-1\right) /b>0$ and with $A_{N}^{*1/N}\rightarrow
a $: $\rho ^{-1}\mathbf{P}\left( K_{N}^{*}=0\right) \sim 1-\frac{1}{b}%
A_{N}^{*}\rightarrow 1.$

\subsubsection{$\theta =0$ \textbf{(infinite mean case)}}

Let $a_{n}\in \left( 0,1\right) $. The case $\theta =0$ is defined by
continuity by $f_{n}\left( z\right) =1-\left( 1-\rho _{n}\right)
^{1-a_{n}}\left( 1-z\right) ^{a_{n}}$, provided $b_{n}$ is given by $\rho
_{n}=1-\left( \left( 1-a_{n}\right) /b_{n}\right) ^{1/\theta }\in \left(
0,1\right) $. Here $f_{n}^{\prime }\left( 1\right) =\infty $. Let $\lambda
_{n}:=\left( 1-\rho _{n}\right) ^{1-a_{n}}$ and consider then $f_{n}\left(
z\right) =1-\lambda _{n}\left( 1-z\right) ^{a_{n}}$. It holds that 
\begin{equation*}
\phi _{N}^{*}\left( z\right) =1-\Lambda _{N}^{*}\left( 1-z\right)
^{A_{N}^{*}},
\end{equation*}
where $\Lambda _{N}^{*}=\prod_{n=1}^{N}\lambda _{n}^{\prod_{m=n+1}^{N}a_{m}}$%
. Thus, $\mathbf{E}\left( K_{N}^{*}\right) =\infty $ (strong
supercriticality of such a coalescence process) and 
\begin{equation*}
\mathbf{P}\left( K_{N}^{*}=0\right) =\phi _{N}^{*}\left( 0\right) =1-\Lambda
_{N}^{*}.
\end{equation*}
We can produce an example for which\textbf{\ }$\mathbf{P}\left(
K_{N}^{*}=0\right) \underset{N\rightarrow \infty }{\rightarrow }\rho >0$%
\textbf{\ }and control its speed of convergence. Take $a_{n}=a\left(
1-\varepsilon _{n}\right) $ with $\varepsilon _{n}\rightarrow 0$ fast enough
($a\in \left( 0,1\right) $) and suppose $\lambda _{n}=\lambda ^{1-a_{n}}$
where $\lambda \in \left( 0,1\right) $. Then 
\begin{equation*}
1-\phi _{N}^{*}\left( z\right) =\lambda \left( \frac{1-z}{\lambda }\right)
^{A_{N}^{*}},
\end{equation*}
showing that (again because $A_{N}^{*1/N}\rightarrow a$) 
\begin{equation}
\mathbf{P}\left( K_{N}^{*}=0\right) =\phi _{N}^{*}\left( 0\right) \sim
1-\lambda ^{1-A_{N}^{*}}\rightarrow 1-\lambda =:\rho >0.  \label{ebr3}
\end{equation}
If for instance $\varepsilon _{n}\rightarrow 0$ like $n^{-\alpha }$ with $%
\alpha >1$, $\mathbf{P}\left( K_{N}^{*}>0\right) \rightarrow \lambda =1-\rho
>0$ at double exponential speed: $\lambda ^{-1}\mathbf{P}\left(
K_{N}^{*}>0\right) \sim \lambda ^{-Ca^{n}}$, for some constant $C>0$.

\subsection{The quadratic coalescing mechanism}

Consider the quadratic map: $f\left( z\right) =a\left( b+z\right) ^{2}-b.$
We have $f\left( z\right) =h^{-1}\left( g\left( h\left( z\right) \right)
\right) $ where $h\left( z\right) =z+b$, $g\left( z\right) =az^{2}.$ For
such maps $f$, 
\begin{equation}
f^{\circ 0}\left( z\right) :=z\text{ and }f^{\circ n}\left( z\right)
=a^{2^{n}-1}\left( b+z\right) ^{2^{n}}-b\text{, }n\geq 1  \label{iter}
\end{equation}
and the $n$-th iterates $f^{\circ n}\left( \cdot \right) =f\left( f\left(
...f\left( \cdot \right) \right) \right) $ ($n$ times) of $f$ can explicitly
be found as a degree $2^{n}$-polynomial. Note that $\left( f^{\circ
n}\right) ^{-1}\left( z\right) =f^{\circ -n}\left( z\right) $.

Can $f\left( z\right) $ be a non trivial (ie different from the purely
quadratic map $az^{2}$) pgf so that $f^{\circ n}\left( z\right) $ would be a
pgf for all $n$? For this to hold, we should have $a,b>0$ and $ab\geq 1$ but
also $f\left( 0\right) \in \left( 0,1\right) $ and $f\left( 1\right) =1.$
The latter conditions require that $a=1/\left( b+1\right) $ which cannot be
fulfilled together with $ab\geq 1$. So if $a,b>0$ and $ab\geq 1$, $f\left(
z\right) $ is a quadratic map with positive $\left[ z^{k}\right] -$%
coefficients, $k=0,1,2$, but these summing to a constant larger than one. So
these coefficients cannot be probabilities and neither such a $f$ nor its $n$%
-th iterate can be a genuine pgf. Observe that the $\left[ z^{k}\right] -$%
coefficients of $f^{\circ n}\left( z\right) $ are also non-negative in this
case for all $n$.

We nevertheless have the following result:

\begin{proposition}
$\phi _{n}^{*}\left( z\right) =f^{\circ n}\left( z\right) /f^{\circ n}\left(
1\right) $, as given from (\ref{iter}), is the explicit pgf of a
time-inhomogeneous coalescing Markov process whose generating coalescing
mechanism is 
\begin{equation*}
f_{n+1}\left( z\right) =\frac{\left( ab+\left( \left( a\left( b+1\right)
\right) ^{2^{n}}-ab\right) z\right) ^{2}-ab}{\left( a\left( b+1\right)
\right) ^{2^{n+1}}-ab}.
\end{equation*}
\end{proposition}

\textbf{Proof:} suppose that $f$ is chosen with $a,b>0$ and $ab\geq 1$ so
that it is absolutely monotone on $\left( 0,1\right) $. Consider the pgf $%
\phi _{n}^{*}\left( z\right) =f^{\circ n}\left( z\right) /f^{\circ n}\left(
1\right) $, normalizing $f^{\circ n}\left( z\right) $ so as $\phi
_{n}^{*}\left( 1\right) =1$, leading to 
\begin{equation}
\phi _{n}^{*}\left( z\right) =\frac{a^{2^{n}-1}\left( b+z\right) ^{2^{n}}-b}{%
a^{2^{n}-1}\left( b+1\right) ^{2^{n}}-b}=\frac{\left[ a\left( b+z\right)
\right] ^{2^{n}}-ab}{\left[ a\left( b+1\right) \right] ^{2^{n}}-ab}.
\label{3.0}
\end{equation}
We have $\phi _{0}^{*}\left( z\right) =z$ and 
\begin{equation}
\phi _{n+1}^{*}\left( z\right) =f^{\circ n+1}\left( z\right) /f^{\circ
n+1}\left( 1\right) =\frac{f\left( f^{\circ n}\left( 1\right) \phi ^{\circ
n}\left( z\right) \right) }{f\left( f^{\circ n}\left( 1\right) \right) }%
=:f_{n+1}\left( \phi _{n}^{*}\left( z\right) \right) ,  \label{3.1}
\end{equation}
where $f_{n+1}\left( z\right) :=\frac{f\left( f^{\circ n}\left( 1\right)
z\right) }{f\left( f^{\circ n}\left( 1\right) \right) }$ is an inhomogeneous
pgf, with $f^{\circ n}\left( 1\right) =a^{2^{n}-1}\left( b+1\right)
^{2^{n}}-b$, $f\left( f^{\circ n}\left( 1\right) \right) =f^{\circ
n+1}\left( 1\right) ,$ explicitly known. The inhomogeneous binary coalescing
mechanism of such a Markov process reads 
\begin{eqnarray*}
f_{n+1}\left( z\right) &=&\frac{f\left( \left( a^{2^{n}-1}\left( b+1\right)
^{2^{n}}-b\right) z\right) }{a^{2^{n+1}-1}\left( b+1\right) ^{2^{n+1}}-b}=%
\frac{a\left( b+\left( a^{2^{n}-1}\left( b+1\right) ^{2^{n}}-b\right)
z\right) ^{2}-b}{a^{2^{n+1}-1}\left( b+1\right) ^{2^{n+1}}-b} \\
&=&\frac{a^{2}\left( b+\left( a^{2^{n}-1}\left( b+1\right) ^{2^{n}}-b\right)
z\right) ^{2}-ab}{a^{2^{n+1}}\left( b+1\right) ^{2^{n+1}}-ab}=\frac{\left(
ab+\left( \left( a\left( b+1\right) \right) ^{2^{n}}-ab\right) z\right)
^{2}-ab}{\left( a\left( b+1\right) \right) ^{2^{n+1}}-ab}.
\end{eqnarray*}

Note $f_{n+1}\left( z\right) \underset{n\rightarrow \infty }{\rightarrow }%
z^{2}$ because $a\left( b+1\right) >1$: the asymptotic inhomogeneous driving
pgf is purely quadratic (a purely binary coalescence process) with $%
f_{n+1}\left( 0\right) =ab\left( ab-1\right) /\left( \left( a\left(
b+1\right) \right) ^{2^{n+1}}-ab\right) \rightarrow 0$, very fast. $\square $

For such a binary coalescing process, we have:\newline

$\left( i\right) $ $\phi _{N}^{*\prime }\left( 1\right) =\mathbf{E}\left(
K_{N}^{*}\right) =\frac{a^{2^{N}-1}2^{N}\left( b+1\right) ^{2^{N}-1}}{%
a^{2^{N}-1}\left( b+1\right) ^{2^{N}}-b}\sim \frac{2^{N}}{b+1}\rightarrow
\infty $. Thus $\mathbf{E}\left( K_{N}^{*}\right) ^{1/N}\rightarrow \mu =2>1$%
. The binary coalescent process is supercritical.\newline

$\left( ii\right) $ Concerning the probability of an empty box, $\mathbf{P}%
\left( K_{N}^{*}=0\right) =\phi _{N}^{*}\left( 0\right) =\left( \left(
ab\right) ^{2^{N}}-ab\right) /\left( \left( a\left( b+1\right) \right)
^{2^{N}}-ab\right) $ $\underset{N\rightarrow \infty }{\sim }\left( b/\left(
b+1\right) \right) ^{2^{N}}\rightarrow \rho =0$. The sequence $\phi
_{N}^{*}\left( 0\right) =\mathbf{P}\left( K_{N}^{*}=0\right) $\textbf{\ }is
decreasing very fast with $N$\ (at double exponential speed) and so $%
K_{n}^{*}=0\nRightarrow K_{n+1}^{*}=0$: the coalescent process can hit $0$\
where it can be regenerated. In addition, for $j\in \left\{
0,...,2^{n}\right\} $, 
\begin{equation}
\mathbf{P}\left( K_{n}^{*}=j\right) =\left[ z^{j}\right] \phi _{n}^{*}\left(
z\right) =\frac{a^{2^{n}-1}}{a^{2^{n}-1}\left( b+1\right) ^{2^{n}}-b}\binom{%
2^{n}}{j}b^{2^{n}-j}  \label{bfl}
\end{equation}
gives the full distribution of $K_{n}^{*}$.\newline

\textbf{Acknowledgments: }TH acknowledges support from the ``Chaire \textit{%
Mod\'{e}lisation math\'{e}matique et biodiversit\'{e}''.} NG and TH also
acknowledge support from the labex MME-DII Center of Excellence (\textit{%
Mod\`{e}les math\'{e}matiques et \'{e}conomiques de la dynamique, de
l'incertitude et des interactions}, ANR-11-LABX-0023-01 project). The
authors are indebted to J. Avan for careful reading of the manuscript.

\end{document}